\documentclass{amsart}

\usepackage{amsthm}
\usepackage[leqno]{amsmath}
\usepackage{latexsym,amsfonts,amssymb}
\usepackage[all]{xy} \SelectTips{eu}{} \SilentMatrices
\usepackage[hyperfootnotes=false]{hyperref}

\newcommand{\numberseries}{\mdseries}   

\newlength{\thmtopspace}                
\newlength{\thmbotspace}                
\newlength{\thmheadspace}               
\newlength{\thmindent}                  

\renewcommand{\subparagraph}{\vspace*{\thmbotspace}}

\newtheoremstyle{bfupright head,slanted body}
                {\thmtopspace}{\thmbotspace}
                {\slshape}{\thmindent}{\bfseries}{.}{\thmheadspace}
                {{\numberseries \thmnumber{(#2) }}\thmnote{#3}}

\newtheoremstyle{bfupright head,upright body}
                {\thmtopspace}{\thmbotspace}
                {\upshape}{\thmindent}{\bfseries}{.}{\thmheadspace}
                {{\numberseries \thmnumber{(#2) }}\thmnote{#3}}

\newtheoremstyle{bfit head,upright body}
                {\thmtopspace}{\thmbotspace}
                {\upshape}{\thmindent}{\upshape}{.}{\thmheadspace}
                {{\numberseries\thmnumber{(#2) }}
                {\bfseries\itshape\thmnote{\negthickspace#3}}}

\newtheoremstyle{it head,upright body}
                {\thmtopspace}{\thmbotspace}
                {\upshape}{\thmindent}{\upshape}{.}{\thmheadspace}
                {{\numberseries\thmnumber{(#2) }}
                {\itshape\thmnote{\negthickspace#3}}}

\newtheoremstyle{fixed bf head,slanted body}
                {\thmtopspace}{\thmbotspace}{\slshape}
                {\thmindent}{\bfseries}{.}{\thmheadspace}
                {{\numberseries \thmnumber{(#2) }}\thmname{#1}\thmnote{ (#3)}}

\newtheoremstyle{fixed bf head,upright body}
                {\thmtopspace}{\thmbotspace}{\upshape}
                {\thmindent}{\bfseries}{.}{\thmheadspace}
                {{\numberseries \thmnumber{(#2) }}\thmname{#1}\thmnote{ (#3)}}

\newtheoremstyle{fixed bfit head,upright body}
                {\thmtopspace}{\thmbotspace}{\upshape}
                {\thmindent}{\bfseries\itshape}{.}{\thmheadspace}
                {{\numberseries \thmnumber{(#2) }}\thmname{#1}\thmnote{ (#3)}}

\newtheoremstyle{sc head,small body}
                {\thmtopspace}{\thmbotspace}
                {\small\upshape}{\thmindent}{\scshape}{.}{\thmheadspace}
                {\thmname{#1}}

\newtheoremstyle{numbered paragraph}
                {\thmtopspace}{\thmbotspace}{\upshape}
                {\thmindent}{\upshape}{}{0pt}
                {{\numberseries \thmnumber{(#2) }}}

\newtheoremstyle{unnumbered paragraph}
                {\thmtopspace}{\thmbotspace}{\upshape}
                {\parindent}{\upshape}{}{0pt}

\theoremstyle{bfupright head,slanted body}
\newtheorem{res}{}[section]             \newtheorem*{res*}{}

\theoremstyle{bfit head,upright body}
                 \newtheorem*{com*}{}

\theoremstyle{bfupright head,upright body}
\newtheorem{bfhpg}[res]{}               \newtheorem*{bfhpg*}{}

\theoremstyle{it head,upright body}
               \newtheorem*{ithpg*}{}

\theoremstyle{sc head,small body}

\theoremstyle{fixed bf head,slanted body}
\newtheorem{thm}[res]{Theorem}          \newtheorem*{thm*}{Theorem}
\newtheorem{prp}[res]{Proposition}      \newtheorem*{prp*}{Proposition}
\newtheorem{cor}[res]{Corollary}        \newtheorem*{cor*}{Corollary}
\newtheorem{lem}[res]{Lemma}            \newtheorem*{lem*}{Lemma}

\theoremstyle{fixed bf head,upright body}
\newtheorem{dfn}[res]{Definition}       \newtheorem*{dfn*}{Definition}
     \newtheorem*{con*}{Construction}
      \newtheorem*{obs*}{Observation}
\newtheorem{rmk}[res]{Remark}           \newtheorem*{rmk*}{Remark}
          \newtheorem*{exa*}{Example}
         \newtheorem*{exe*}{Exercise}
            \newtheorem{stp*}{Setup}

\theoremstyle{numbered paragraph}
\newtheorem{ipg}[res]{}

\theoremstyle{unnumbered paragraph}
\newtheorem{ipg*}{}

\newlength{\thmlistleft}        
\newlength{\thmlistright}       
\newlength{\thmlistpartopsep}   
\newlength{\thmlisttopsep}      
\newlength{\thmlistparsep}      
\newlength{\thmlistitemsep}     

\newcounter{eqc} 
\newenvironment{eqc}{\begin{list}{\upshape (\textit{\roman{eqc}})}%
    {\usecounter{eqc}%
      \setlength{\leftmargin}{\thmlistleft}%
      \setlength{\labelwidth}{\thmlistleft}%
      \setlength{\rightmargin}{\thmlistright}%
      \setlength{\partopsep}{\thmlistpartopsep}%
      \setlength{\topsep}{\thmlisttopsep}%
      \setlength{\parsep}{\thmlistparsep}%
      \setlength{\itemsep}{\thmlistitemsep}}}%
  {\end{list}}%

\newcounter{prt}
  {\end{list}}%

\newcounter{rqm}
\newenvironment{rqm}{\begin{list}{\upshape (\arabic{rqm})}%
    {\usecounter{rqm}%
      \setlength{\leftmargin}{\thmlistleft}%
      \setlength{\labelwidth}{\thmlistleft}%
      \setlength{\rightmargin}{\thmlistright}%
      \setlength{\partopsep}{\thmlistpartopsep}%
      \setlength{\topsep}{\thmlisttopsep}%
      \setlength{\parsep}{\thmlistparsep}%
      \setlength{\itemsep}{\thmlistitemsep}}}%
  {\end{list}}%

  {\end{list}}%

  {\end{list}}%

\newenvironment{prf*}[1][Proof]{%
  \begin{proof}[\bf #1]
    \setcounter{equation}{0}
    \renewcommand{\theequation}{\arabic{equation}}}
  {\end{proof}
}

\newcommand{\pgref}[1]{(\ref{#1})}

\newcommand{\thmref}[2][Theorem~]{#1\pgref{thm:#2}}
\newcommand{\corref}[2][Corollary~]{#1\pgref{cor:#2}}
\newcommand{\prpref}[2][Proposition~]{#1\pgref{prp:#2}}
\newcommand{\lemref}[2][Lemma~]{#1\pgref{lem:#2}}

\newcommand{\dfnref}[2][Definition~]{#1\pgref{dfn:#2}}

\renewcommand{\eqref}[1]{\pgref{eq:#1}}

\newcommand{\thmcite}[2][?]{\cite[Thm.~#1]{#2}}
\newcommand{\corcite}[2][?]{\cite[Cor.~#1]{#2}}
\newcommand{\prpcite}[2][?]{\cite[Prop.~#1]{#2}}

\newcommand{\set}[1]{\{\,#1\,\}}
\newcommand{\setof}[2]{\{\,#1 \mid #2\,\}}

\newcommand{\ZZ}{\mathbb{Z}}

\newcommand{\qtext}[1]{\quad\text{#1}\quad}
\newcommand{\qand}{\qtext{and}}
\newcommand{\deq}{\:=\:}

\renewcommand{\a}{\alpha}
\renewcommand{\b}{\beta}
\newcommand{\f}{\varphi}
\newcommand{\g}{\gamma}
\newcommand{\m}{\mathfrak{m}}
\newcommand{\n}{\mathfrak{n}}
\newcommand{\q}{\mathfrak{q}}
\newcommand{\z}{\zeta}

\newcommand{\ot}{\otimes}

\newcommand{\is}{\cong}

\renewcommand{\le}{\leqslant}
\renewcommand{\ge}{\geqslant}

\newcommand{\into}{\hookrightarrow}

\newcommand{\lra}{\longrightarrow}

\newcommand{\xra}[2][]{\xrightarrow[#1]{\;#2\;}}

\newcommand{\ira}{\xra{\;\is\;}}

\newcommand{\dra}[2][]{\xra{\dif[#1]{#2}}}

\newcommand{\pows}[2][k]{#1[\mspace{-2.3mu}[#2]\mspace{-2.3mu}]}

\newcommand{\Rmk}{(R,\m,k)}

\newcommand{\Snl}{(S,\n,l)}

\newcommand{\Rhat}{\widehat{R}}

\newcommand{\mapdef}[4][\rightarrow]{\nobreak{#2\colon #3 #1 #4}}
\newcommand{\dmapdef}[4][\lra]{\nobreak{#2\colon #3\:#1\:#4}}
\renewcommand{\Im}[1]{\nobreak{\operatorname{Im}#1}}
\newcommand{\Ker}[1]{\nobreak{\operatorname{Ker}#1}}
\newcommand{\Coker}[1]{\nobreak{\operatorname{Coker}#1}}

\newcommand{\dif}[2][]{{\partial}^{#2}_{#1}}

\newcommand{\dptR}{\operatorname{depth}R}
\newcommand{\dimR}{\operatorname{dim}R}

\newcommand{\Hom}[3][R]{\operatorname{Hom}_{#1}(#2,#3)}
\newcommand{\Ext}[4][R]{\operatorname{Ext}_{#1}^{#2}(#3,#4)}
\newcommand{\tp}[3][R]{\nobreak{#2\otimes_{#1}#3}}
\newcommand{\tpp}[3][R]{\left(\tp[#1]{#2}{#3}\right)}

\newcommand{\Tor}[4][R]{\operatorname{Tor}^{#1}_{#2}(#3,#4)}

\makeatletter
\def\@nobreak@#1{\mathchoice%
  {\nobreakdef@\displaystyle\f@size{#1}}%
  {\nobreakdef@\nobreakstyle\tf@size{\firstchoice@false #1}}%
  {\nobreakdef@\nobreakstyle\sf@size{\firstchoice@false #1}}%
  {\nobreakdef@\nobreakstyle\ssf@size{\firstchoice@false #1}}%
  \check@mathfonts}%
\def\nobreakdef@#1#2#3{\hbox{{%
                    \everymath{#1}%
                    \let\f@size#2\selectfont%
                    #3}}}%
\makeatother

\numberwithin{equation}{res}

\renewcommand{\Rmk}{(R,\m,\k)}
\renewcommand{\Snl}{(S,\n,\l)}

\newcommand{\be}{\varepsilon}
\newcommand{\dd}{\partial}
\newcommand{\ddK}{\mathrm{d}}
\newcommand{\Kc}{\mathrm{K}}
\newcommand{\G}[1][R]{\mathcal{G}(#1)}
\newcommand{\F}[1][R]{\mathcal{F}(#1)}

\newcommand{\Gbar}[1][R]{\tp[]{S}{\mathcal{G}(#1)}}
\renewcommand{\bar}[2][S]{\tp[]{#1}{\cat{#2}}}

\newcommand{\cl}[1]{\langle#1\rangle}
\newcommand{\add}[1]{\operatorname{add}(#1)}
\newcommand{\Gbarcl}[1][R]{\cl{\Gbar}}

\renewcommand{\k}{\mathsf{k}}
\renewcommand{\l}{\ell}

\newcommand{\syz}[2]{{#2}_{#1}}
\newcommand{\syzp}[2]{({#2}_{#1})}
\newcommand{\Syz}[3][R]{\Omega^{#1}_{#2}(#3)}
\newcommand{\x}{\pmb{x}}
\newcommand{\y}{\pmb{y}}
\newcommand{\cat}[1]{\mathcal{#1}}
\renewcommand{\mod}[1]{\operatorname{mod}(#1)}

\newcommand{\one}{\ensuremath{\mathord \ast}}

\renewcommand{\H}[2][]{\operatorname{H}^{#1}(#2)}

\begin{document}
\enlargethispage{\baselineskip}
\vspace*{-.1\baselineskip}

\title[Finite Gorenstein representation type implies simple
singularity]{Finite Gorenstein representation type\\implies simple
  singularity}

\dedicatory{To Lucho Avramov on his sixtieth birthday}

\author[L.~W.~Christensen]{Lars Winther Christensen}

\address{Department of Math.\ and Stat., Texas Tech University,
  Lubbock, TX~79409, U.S.A.}

\email{lars.w.christensen@ttu.edu}%

\thanks{This work was done while L.W.C.\ visited University of
  Nebraska--Lincoln (UNL), partly supported by a grant from the
  Carlsberg Foundation. Part of it was done during R.T.'s visit to UNL
  supported by NSF grant DMS~0201904. J.S.\ was supported by NSF grant
  DMS~0201904\vspace{1ex}}%

\author[G. Piepmeyer]{Greg Piepmeyer}

\address{Department of Mathematics, University of Nebraska, Lincoln,
  NE~68588, U.S.A.}

\curraddr{Department of Math., University of Missouri, Columbia,
  MO~65211, U.S.A.}

\email{greg@math.missouri.edu}

\author[J. Striuli]{Janet Striuli}

\address{Department of Mathematics, University of Nebraska, Lincoln,
  NE~68588, U.S.A.}

\email{jstriuli2@math.unl.edu}

\author[R.~Takahashi]{Ryo~Takahashi}

\address{Department of Mathematical Sciences, Faculty of Science,
  Shinshu University 3-1-1 Asahi, Matsumoto, Nagano 390-8621, Japan}

\email{takahasi@math.shinshu-u.ac.jp}

\date{22 February 2008}




\begin{abstract}
  Let $R$ be a commutative noetherian local ring and consider the set
  of isomorphism classes of indecomposable totally reflexive
  $R$-modules. We prove that if this set is finite, then either it has
  exactly one element, represented by the rank $1$ free module, or $R$
  is Gorenstein and an isolated singularity (if $R$ is complete, then
  it is even a simple hypersurface singularity).  The crux of our
  proof is to argue that if the residue field has a totally reflexive
  cover, then $R$ is Gorenstein or every totally reflexive $R$-module
  is free.
\end{abstract}

\maketitle

\thispagestyle{empty}

\vspace{-.5\baselineskip}
\section*{Introduction}

Remarkable connections between the module theory of a local ring and
the character of its singularity emerged in the 1980s. They show how
finiteness conditions on the category of maximal Cohen--Macaulay
modules\footnote{ The finitely generated modules whose depth equals
  the Krull dimension of the ring.}  characterize particular isolated
singularities. We develop these connections in several directions.

A local ring with only finitely many isomorphism classes of
indecomposable maximal Cohen--Macaulay modules is said to be of finite
Cohen--Macaulay (CM) representation type.  By work of Auslander
\cite{MAs85}, every complete Cohen--Macaulay local ring of finite
CM~representation type is an isolated singularity.

Specialization to Gorenstein rings opens to a finer description of the
\mbox{singularities}; it centers on the simple hypersurface
singularities identified in Arnol$'$d's work on germs of holomorphic
functions \cite{VIA75}.  By work of Buchweitz, Greuel, and Schreyer
\cite{BGS-87}, Herzog \cite{JHr78}, and Yoshino \cite{yos}, a complete
Gorenstein ring of finite CM represen\-tation type is a simple
singularity in the generalized sense of \cite{yos}.  Under extra
assumptions on the ring, the converse holds by work of Kn\"{o}rrer
\cite{HKn87} and Solberg~\cite{OSl89}.

In this introduction, $R$ is a commutative noetherian local ring with
maximal ideal~$\m$ and residue field $\k$.  To avoid the \emph{a
  priori} condition in \cite{BGS-87,JHr78,yos} that $R$ is Gorenstein,
we replace finite CM repre\-sentation type with a finiteness condition
on the category $\G$ of modules of Gorenstein dimension~$0$.  Over a
Gorenstein ring, these modules are precisely the maximal
Cohen--Macaulay modules, but they are known to exist over any ring,
unlike maximal Cohen--Macaulay modules.
\begin{res*}[Theorem~A]
  Let $R$ be complete. If the set of isomorphism classes of
  \mbox{non-free} indecomposable modules in $\G$ is finite and not
  empty, then $R$ is a simple \mbox{singularity}.
\end{res*}
The category $\G$ was introduced by Auslander and Bridger
\cite{MAs67,MAsMBr69}. An $R$-module $G$ is in $\G$ if there is an
exact complex of finitely generated free $R$-modules
\begin{equation*}
  \pmb{F} \deq \cdots \lra F_{n+1} \dra[n+1]{} F_{n} \dra[n]{} F_{n-1}
  \lra \cdots,
\end{equation*}
such that $G$ is isomorphic to $\Coker{\dif[0]{}}$ and the complex
$\Hom{\pmb{F}}{R}$ is exact. Every finitely generated free $R$-module
is in $\G$, and the modules in this category have Gorenstein dimension
$0$ as in \cite{MAs67,MAsMBr69}; following \cite{LLAAMr02} we call
them \emph{totally reflexive}.

The aforementioned works \cite{BGS-87, JHr78, yos} show that Theorem~A
follows from the next result, which is proved as \thmref[]{B}.

\begin{res*}[Theorem~B]
  If the set of isomorphism classes of indecomposable modules in $\G$
  is finite, then $R$ is Gorenstein or every module in $\G$ is free.
\end{res*}

As this theorem does not require $R$ to be complete, we considerably
strengthen Theorem~A using work of Huneke, Leuschke, and R.~Wiegand
\cite{CHnGJL02, GLsRWg00, RWg98}; this occurs in \thmref[]{A}.
Theorem~B was conjectured by R.~Takahashi \cite{RTk05a}, who proved it
for henselian rings of depth at most two
\cite{RTK04a,RTK04d,RTk05a}. The class of rings over which all totally
reflexive modules are free is poorly understood, but it is known to
include all Golod rings \cite{LLAAMr02}, in particular, all
Cohen--Macaulay rings of minimal~multiplicity.

\begin{ipg*}
  To prove Theorem~B we use a notion of \emph{$\G$-approximations},
  which is close kin to the CM-approximations of Auslander and
  Buchweitz \cite{MAsROB89}. When $R$ is Gorenstein, a
  $\G$-approximation is exactly a CM-approximation. By
  \cite{MAsROB89}, every module over a Gorenstein ring has a
  CM-approximation. Our proof of Theorem~B goes via the following
  strong converse, proved as~\thmref[]{B-approximation}.
\end{ipg*}

\begin{res*}[Theorem~C]
  Let $R$ be a local ring and assume there is a non-free module in
  $\G$. If the residue field $\k$ has a $\G$-approximation, then $R$
  is Gorenstein.
\end{res*}

This theorem complements recent developments in relative homological
algebra. The notion of totally reflexive modules has two extensions to
non-finitely generated modules; see \cite{lnm} for details.  One is
Gorenstein projective modules, which allows arbitrary free modules in
the definition above.  By recent work of J{\o}rgensen \cite{PJr07},
every module over a complete local ring has a Gorenstein projective
precover.  The other extension is Gorenstein flat modules.  By a
result of Enochs and L\'{o}pez-Ramos \cite{EEnJLR02}, every module has
a Gorenstein flat precover.

Theorem~C counterposes these developments; it shows that for finitely
generated modules, the precovers found in \cite{PJr07} and
\cite{EEnJLR02} cannot, in general, be finitely generated.  Assume
that $R$ is complete. Then a finitely generated $R$-module has a
$\G$-approximation if and only if it has a $\G$-precover. Assume
further that $R$ is not Gorenstein.  Theorem~C shows that if $X \to
\k$ is a Gorenstein projective/flat precover and $X$ is not free, then
$X$ is not finitely generated.

\enlargethispage*{1.5\baselineskip}

\section{Categories and covers}
\label{sec:1}

In this paper, rings are commutative and noetherian; modules are
finitely generated (unless otherwise specified). We write $\mod{R}$
for the category of finitely generated modules over a ring $R$.

For an $R$-module $M$, we denote by $M_i$ the $i$th syzygy in a free
resolution. When $R$ is local, we denote by $\Syz{i}{M}$ the $i$th
syzygy in the minimal free resolution of~$M$.  For an $R$-module $M$,
set $M^* =\Hom{M}{R}$; we refer to this module as the \emph{algebraic
  dual} of $M$.

We only consider full subcategories of $\mod{R}$; this allows us to
define a subcategory by specifying its objects.  In the following,
$\cat{B}$ is a subcategory of $\mod{R}$.

\begin{bfhpg}[Closures]
  Recall that the category $\cat{B}$ is said to be closed under
  extensions if for every short exact sequence $0 \to B \to X \to B'
  \to 0$ with $B$ and $B'$ in $\cat{B}$ also $X$ is in $\cat{B}$.  The
  closure of $\cat{B}$ under extensions is by definition the smallest
  subcategory containing $\cat{B}$ and closed under extensions.
  Recall also that $\cat{B}$ is closed under direct sums and direct
  summands when a direct sum $M\oplus N$ is in $\cat{B}$ if and only
  if both summands are in $\cat{B}$. The closure of $\cat{B}$ under
  addition is by definition the smallest subcategory containing
  $\cat{B}$ and closed under direct sums and direct summands; we
  denote it by $\add{\cat{B}}$.
  
  We define the closure $\cl{\cat{B}}$ to be the smallest subcategory
  containing $\cat{B}$ and closed under direct summands and
  extensions.  It is straightforward to verify that the closure
  $\cl{\cat{B}}$ is reached by countable alternating iteration,
  starting with $\cat{B}$, between closure under addition and closure
  under extensions.
  
  We say that $\cat{B}$ is \emph{closed under algebraic duality} if
  for every module $B$ in $\cat{B}$ the module $B^*$ is also in
  $\cat{B}$.  Similarly, we say that $\cat{B}$ is \emph{closed under
    syzygies} if for every module $B$ in $\cat{B}$ every first syzygy
  $B_1$ is in $\cat{B}$; then every syzygy $B_i$ is in $\cat{B}$.
\end{bfhpg}

\begin{bfhpg}[Precovers and covers]
  Let $M$ be an $R$-module. A \emph{$\cat{B}$-precover} of $M$ is a
  homomorphism $\f \colon B \to M$, with $B \in \cat{B}$, such that
  every homomorphism $X \to M$ with $X \in \cat{B}$, factors through
  $\f$; i.e., the homomorphism
  \begin{equation*}
    \dmapdef{\Hom{X}{\f}}{\Hom{X}{B}}{\Hom{X}{M}}
  \end{equation*}
  is surjective for each module $X$ in $\cat{B}$. A $\cat{B}$-precover
  $\f \colon B \to M$ is a \emph{$\cat{B}$-cover} if every $\g \in
  \Hom{B}{B}$ with $\f\g=\f$ is an automorphism.
  
  Note that if the category $\cat{B}$ contains $R$, then every
  $\cat{B}$-precover is surjective.
\end{bfhpg}

\begin{ipg}
  \label{pc}
  If there are only finitely many isomorphism classes of
  indecomposable modules in $\cat{B}$, then every finitely generated
  $R$-module has a $\cat{B}$-precover; see~\prpcite[4.2]{MAsSOS80}.
\end{ipg}

\begin{ipg}
  \label{summand}
  Consider a diagram $\smash{\xymatrix@C=2.25em{B \ar[r]^-{\f} & M
      \oplus N \ar@<.5ex>[r]^-{\pi} & M \ar@<.5ex>[l]^-{\iota}}}$,
  where $\pi\iota$ is the identity on $M$. If $\f$ is a
  $\cat{B}$-precover, then so is $\mapdef{\pi\f}{B}{M}$.
\end{ipg}

The next two lemmas appear in Xu's book \cite[2.1.1 and 1.2.8]{xu}.
We include a proof of the second one since Xu left it to the reader.

\begin{bfhpg}[Wakamatsu's lemma]
  \label{lem:wak}
  Let $\cat{B}$ be a subcategory of $\mod{R}$, and let $\f$ be a
  $\cat{B}$-cover of an $R$-module $M$.  If $\cat{B}$ is closed under
  extensions, then $\Ext{1}{X}{\ker{\f}}=0$ for all $X\in \cat{B}$.
\end{bfhpg}

\begin{lem}
  \label{lem:xu}
  Let $\cat{B}$ be a subcategory of $\,\mod{R}$, and let $M$ be an
  $R$-module. If $M$ has a $\cat{B}$-cover, then a $\cat{B}$-precover
  $\f \colon X \to M$ is a cover if and only if $\,\ker \f$ contains
  no non-zero direct summand of $X$.
\end{lem}

\begin{prf*}
  Let $\mapdef{\psi}{Y}{M}$ be a $\cat{B}$-cover. For the ``if'' part,
  consider the commutative diagram below, where $\a$ and $\b$ are
  given by the precovering properties of $\f$ and $\psi$.
  \begin{equation*}
    \xymatrix{
      & M
      \\
      Y
      \ar[r]^-{\a}
      \ar[ru]^-{\psi}
      & X
      \ar[r]^-{\b}
      \ar[u]^-{\f}
      & Y
      \ar[lu]_-{\psi}
    }
  \end{equation*}
  Since $\psi\b\a = \psi$ and $\psi$ is a cover, the composite $\b\a$
  is an automorphism, so $\b$ is surjective. It also follows that $X$
  is isomorphic to $\Ker{\b} \oplus \Im{\a}$. As $\Ker{\f}$ contains
  no non-zero summand of $X$, the inclusion $\Ker{\b} \subseteq
  \Ker{\f}$ implies that $\b$ is also injective.  Consequently, $\f$
  is a $\cat{B}$-cover.
  
  For the ``only if'' part, consider a decomposition $X=Y\oplus Z$,
  and assume there is an inclusion $Z \subseteq \Ker{\f}$. Let $\pi$
  be the endomorphism of $X$ projecting onto $Y$, then $\f\pi=\f$.
  Since $\f$ is a cover, $\pi$ is an automorphism, whence $Z=0$.
\end{prf*}

\section{Approximations and reflexive subcategories}
\label{sec:section 2}

Stability of (pre-)covers under base change is delicate to track. To
avoid this task, we develop a notion between precover and cover.  The
next definition is in line with that of CM-approximations
\cite{MAsROB89}; for $\G$ it broadens the notion used in
\cite{LLAAMr02}.

\begin{bfhpg}[Definitions]
  Let $\cat{B}$ be a subcategory of $\mod{R}$ and set
  \begin{equation*}
    \cat{B}^\perp = \setof{L\in\mod{R}}{\Ext{i}{B}{L}=0 \text{ for all
        $B\in\cat{B}$ and all $i>0$}}.
  \end{equation*}
  
  Let $M$ be an $R$-module. A \emph{$\cat{B}$-approximation} of $M$ is
  a short exact sequence
  \begin{equation*}
    0 \lra L \lra B \lra M \lra 0,
  \end{equation*}
  where $B$ is in $\cat{B}$ and $L$ is in $\cat{B}^\perp$.
\end{bfhpg}

\begin{ipg}
  \label{cov-approx}
  Let $\cat{B}$ be a subcategory of $\mod{R}$ and $M$ be an
  $R$-module.
  
  (a) If $0 \lra \Ker{\f} \lra B \xra{\f} M \lra 0$ is a
  $\cat{B}$-approximation of $M$, then $\f$ is a \emph{special}
  $\cat{B}$-precover of $M$; see \prpcite[2.1.3]{xu}.
  
  (b) If $B \xra{\f} M$ is a surjective $\cat{B}$-cover, and $\cat{B}$
  is closed under syzygies and extensions, then the sequence $0 \lra
  \Ker{\f} \lra B \xra{\f} M \lra 0$ is a $\cat{B}$-approximation of
  $M$ by Wakamatsu's lemma.
  
  (c) Assume $\mod{R}$ has the Krull--Schmidt property (e.g., $R$ is
  henselian) and $\cat{B}$ is closed under direct summands. The module
  $M$ has a $\cat{B}$-cover if and only if it has a
  $\cat{B}$-precover; see \corcite[2.5]{RTk05a}.
\end{ipg}

The next two results study the behavior of approximations under base
change.

Let $\mapdef{\vartheta}{R}{S}$ be a ring homomorphism. We say that
$\vartheta$ is of finite flat dimension if $S$, viewed as an
$R$-module through $\vartheta$, has a bounded resolution by flat
$R$-modules. We write $\Tor{i>0}{S}{\cat{B}} =0$ if for all
$B\in\cat{B}$, and for all $i>0$, the modules $\Tor{i}{S}{B}$ vanish.
We denote by $\bar{B}$ the subcategory of $S$-modules $\tp{S}{B}$ with
$B\in \cat{B}$.

\begin{lem}
  \label{lem:Gbarclperp}
  Let $R\to S$ be a ring homomorphism of finite flat dimension. Let
  $\cat{B}$ be a subcategory of $\,\mod{R}$ such that
  $\Tor{i>0}{S}{\cat{B}} =0$. If $L \in \cat{B}^\perp$ and
  $\Tor{i>0}{S}{L} =0$, then for every $m \in \ZZ$ and every
  $B\in\cat{B}$ there is an isomorphism
  \begin{equation*}
    \Ext[S]{m}{\tp{S}{B}}{\tp{S}{L}} \is \Tor{-m}{S}{\Hom{B}{L}}.
  \end{equation*}
  In particular, there are isomorphisms $\Hom[S]{\tp{S}{B}}{\tp{S}{L}}
  \is \tp{S}{\Hom{B}{L}}$, and $\tp{S}{L}$ is in
  $\cl{\tp[]{S}{\cat{B}}}^\perp$.
\end{lem}

\begin{prf*}
  Fix $B \in \cat{B}$. Take a free resolution $\pmb{E} \to B$ and a
  bounded flat resolution $\pmb{F} \to S$ over $R$. By the vanishing
  of (co)homology, the induced morphisms
  \begin{equation*}
    \tp{S}{\pmb{E}} \to \tp{S}{B}, \qquad \tp{\pmb{F}}{L} \to \tp{S}{L},
    \qand \Hom{B}{L} \to \Hom{\pmb{E}}{L}
  \end{equation*}
  are homology isomorphisms. In particular, the first one is a free
  resolution of the $S$-module $\tp{S}{B}$. The functors
  $\Hom{\pmb{E}}{-}$ and $\tp{\pmb{F}}{-}$ preserves homology
  isomorphisms. This explains the first, third, and fifth isomorphisms
  below.
  \begin{align*}
    \Ext[S]{m}{\tp{S}{B}}{\tp{S}{L}} 
    &\is \H[m]{\Hom[S]{\tp{S}{\pmb{E}}}{\tp{S}{L}}} \\
    &\is \H[m]{\Hom{\pmb{E}}{\tp{S}{L}}} \\
    &\is \H[m]{\Hom{\pmb{E}}{\tp{\pmb{F}}{L}}} \\
    &\is \H[m]{\tp{\pmb{F}}{\Hom{\pmb{E}}{L}}} \\
    &\is \H[m]{\tp{\pmb{F}}{\Hom{B}{L}}} \\
    &\is \Tor{-m}{S}{\Hom{B}{L}}
  \end{align*}
  The second isomorphism follows from Hom-tensor adjointness, and the
  fourth is tensor evaluation; see \prpcite[II.5.14]{rad}. For $m=0$
  the composite isomorphism reads $\Hom[S]{\tp{S}{B}}{\tp{S}{L}} \is
  \tp{S}{\Hom{B}{L}}$. That $\tp{S}{L}$ is in
  $\cl{\tp[]{S}{\cat{B}}}^\perp$ follows as $\operatorname{Tor}^R_i$
  is zero for $i<0$.
\end{prf*}

\begin{prp}
  \label{prp:C-cover}
  Let $R \to S$ be a ring homomorphism and $\cat{B}$ be a subcategory
  of $\,\mod{R}$. Let $M$ be an $R$-module with a
  $\cat{B}$-approximation $0 \to L \to B \to M \to 0$. If
  $\Tor{i>0}{S}{\cat{B}} =0$ and $\Tor{i>0}{S}{M} =0$, then
  \begin{equation*}
    0 \lra \tp{S}{L} \lra \tp{S}{B} \lra \tp{S}{M} \lra 0
  \end{equation*}
  is an $\cl{\bar{B}}$-approximation.
\end{prp}

\begin{prf*}
  By the assumptions on $\cat{B}$ and $M$, application of the functor
  $\tp{S}{-}$ to the $\cat{B}$-approximation of $M$ yields the desired
  short exact sequence and also equalities $\Tor{i>0}{S}{L}=0$. Now
  \lemref{Gbarclperp} gives that $\tp{S}{L}$ is in
  $\cl{\bar{B}}^\perp$.
\end{prf*}

\begin{ipg}
  Let $\cat{B}$ be a subcategory of $\mod{R}$ with $R \in
  \cat{B}^\perp$. For every $B\in\cat{B}$ and every $R$-module $N$,
  dimension shifting yields
  \begin{equation*}
    \Ext{i}{B}{N} \is \Ext{i+h}{B}{\syz{h}{N}}\quad \text{for
      $i>0$ and $h \ge 0$.}
  \end{equation*}
  Moreover, for $h \ge 0$ the algebraic dual $B^*$ is a $h$th syzygy
  of $(B_h)^*$, so
  \begin{equation*}
    \Ext{i}{B^*}{N} \is \Ext{i+h}{\syzp{h}{B}^*}{N}\quad \text{for
      $i>0$ and $h \ge 0$.}
  \end{equation*}
  If, furthermore, $\cat{B}$ is closed under syzygies and algebraic
  duality, then these isomorphisms combine to yield
  \begin{equation}
    \label{eq:dim shift}
    \Ext{i}{B^*}{\syz{j}{N}} \is 
    \Ext{i}{\syzp{h}{B}^*}{\syz{j-h}{N}}\quad \text{for $i>0$ and $j
      \ge h \ge 0$.}
  \end{equation}
\end{ipg}

In particular, \eqref{dim shift} holds when $\cat{B}$ is a category
satisfying the next definition.

\begin{dfn}
  \label{dfn:reflexive}
  A subcategory $\cat{B}$ of $\mod{R}$ is \emph{reflexive} if $R$ is
  in $\cat{B} \cap \cat{B}^\perp$ and $\cat{B}$ is closed under
  \begin{rqm}
  \item direct sums and direct summands,
  \item syzygies, and
  \item algebraic duality.
  \end{rqm}
\end{dfn}

It is standard that the category $\G$ of totally reflexive $R$-modules
is a reflexive subcategory of $\mod{R}$. Moreover, using the
characterization of $\G$ provided by \cite[(1.1.2) and (4.1.4)]{lnm},
it is straightforward to verify that every reflexive subcategory of
$\mod{R}$ is, in fact, a subcategory of $\G$.

\begin{ipg}
  \label{FBG}
  In the rest of the paper, $\F$ denotes the category of finitely
  generated free \mbox{$R$-modules}. Let $\cat{B}$ be a reflexive
  subcategory of $\mod{R}$. There are containments
  \begin{equation*}
    \F \subseteq \cat{B} \subseteq \G.
  \end{equation*}
  Further, let $R \to S$ be a ring homomorphism of finite flat
  dimension, then
  \begin{equation*}
    \Tor{i>0}{S}{\cat{B}}=0,
  \end{equation*}
  as every module in $\cat{B}$ is an infinite syzygy.
\end{ipg}

The next observation is crucial for our proofs of the main theorems.

\begin{ipg}
  \label{precov-approx}
  Assume $\mod{R}$ has the Krull--Schmidt property (e.g., $R$ is
  henselian) and let $\cat{B}$ be a reflexive subcategory of $\mod{R}$
  closed under extensions. We claim that an $R$-module $M$ has a
  $\cat{B}$-precover if and only if it has a $\cat{B}$-approximation.
  Indeed, let $\mapdef{\f}{B}{M}$ be a $\cat{B}$-precover; by
  \pgref{cov-approx}(c) the module $M$ also has a $\cat{B}$-cover.
  Decompose $B$ as $B' \oplus B''$, where $B''$ is the largest direct
  summand of $B$ contained in $\Ker{\f}$. By \lemref{xu} the
  factorization $\mapdef{\f'}{B'}{M}$ is a cover, and by
  \pgref{cov-approx}(b) the sequence $0 \to \Ker{\f'} \to B' \to M \to
  0$ is a $\cat{B}$-approximation.
\end{ipg}

\begin{lem}
  \label{lem:syz-approximation}
  Let $\cat{B}$ be a reflexive subcategory of $\,\mod{R}$ and $M$ be
  an $R$-module. If $M$ has a $\cat{B}$-approximation, then every
  syzygy of $M$ has a $\cat{B}$-approximation.
\end{lem}

\begin{prf*}
  Let $0 \to L \to B \to M \to 0$ be a $\cat{B}$-approximation. It is
  sufficient to prove that every first syzygy $M_1$ has a
  $\cat{B}$-approximation. By the horseshoe construction, there is a
  short exact sequence $0 \to L_1 \to B_1 \to M_1 \to 0$, and the
  syzygy $B_1$ is in $\cat{B}$ by assumption.  Let $X$ be in
  $\cat{B}$. Since $\cat{B}$ is reflexive, there is an isomorphism $X
  \is X^{**}$, and also the module $((X^*)_1)^*$ is in $\cat{B}$.  Now
  \eqref{dim shift} yields the second isomorphism in the chain
  \begin{equation*}
    \Ext{i}{X}{L_1} \is \Ext{i}{X^{**}}{L_1} \is \Ext{i}{\,((X^*)_1)^*}{L\,}
    = 0.\qedhere
  \end{equation*}
\end{prf*}

\begin{prp}
  \label{prp:laundry list}
  Let $R \to S$ be a ring homomorphism of finite flat dimension. If
  $\cat{B}$ is a reflexive subcategory of $\,\mod{R}$, then
  $\cl{\bar{B}}$ is a reflexive subcategory of $\,\mod{S}$. In
  particular, $\cl{\tp[]{S}{\G}}$ is reflexive.
\end{prp}

\begin{prf*}
  The ring $S$ is in $\cl{\bar{B}}$.  As $R\in\cat{B}^\perp$, it
  follows from \pgref{FBG} and \lemref{Gbarclperp} that $S$ is in
  $\cl{\bar{B}}^\perp$. By definition, $\cl{\bar{B}}$ is closed under
  direct sums and direct summands; this leaves (2) and (3) in
  \dfnref{reflexive} to verify.
  
  First we prove closure under syzygies. Take $B\in\cat{B}$ and
  consider a short exact sequence $0 \to B_1 \to F \to B \to 0$, where
  $F$ is a free $R$-module. By assumption, the syzygy $B_1$ is in
  $\cat{B}$.  By \pgref{FBG} the sequence
  \begin{equation*}
    0 \lra \tp{S}{B_1} \lra \tp{S}{F} \lra \tp{S}{B} \lra 0
  \end{equation*}
  is exact. It shows that the syzygy $\tp{S}{B_1}$ of $\tp{S}{B}$ is
  in $\bar{B}$. Moreover, it follows that any summand of $\tp{S}{B}$
  has a first syzygy in $\add{\bar{B}}$, in particular, in
  $\cl{\bar{B}}$.  By Schanuel's lemma, a module in $\cl{\bar{B}}$
  with some first syzygy in $\cl{\bar{B}}$ has every first syzygy in
  $\cl{\bar{B}}$.  Finally, given a short exact sequence $0 \to M \to
  X \to N \to 0$, where $M$, $N$, and their first syzygies are in
  $\cl{\bar{B}}$, we claim that also a first syzygy of $X$ is in
  $\cl{\bar{B}}$. Indeed, take presentations of $M$ and $N$. Since
  $\cl{\bar{B}}$ is closed under extensions, it follows from the
  horseshoe construction that a first syzygy of $X$ is in
  $\cl{\bar{B}}$.
  
  Next we prove closure under algebraic duality. Take $B \in \cat{B}$
  and note that by \pgref{FBG}, \lemref{Gbarclperp} applies (with
  $L=R$) to yield the isomorphism
  \begin{equation*}
    \Hom[S]{\tp{S}{B}}{S} \is \tp{S}{\Hom{B}{R}}.
  \end{equation*}
  Thus, the algebraic dual of $\tp{S}{B}$ is in $\bar{B}$. Moreover,
  the algebraic dual of any summand of $\tp{S}{B}$ is in
  $\add{\bar{B}}$, in particular, in $\cl{\bar{B}}$. It is now
  sufficient to prove that for every short exact sequence $0 \to M \to
  X \to N \to 0$, where $M$, $N$, and the duals $M^*$ and $N^*$ are in
  $\cl{\bar{B}}$, also the dual $X^*$ is in $\cl{\bar{B}}$.  Since
  $\cl{\bar{B}}$ is closed under extensions, this is immediate from
  the exact sequence
  \begin{equation*}
    0 \to N^* \to X^* \to M^* \to \Ext[S]{1}{N}{S},
  \end{equation*}
  where $\Ext[S]{1}{N}{S}=0$ as $S$ is in $\cl{\bar{B}}^\perp$.
\end{prf*}

\section{Approximations detect the Gorenstein property}
\label{sec:G-results}

The main result of this section is Theorem~C from the introduction.
\lemref{shopping list} furnishes the base case; for that we study a
standard homomorphism.

\begin{ipg}
  \label{nat}
  For modules $X$ and $N$ over a ring $S$ there is a natural map
  \begin{equation*}
    \dmapdef{\theta_{XN}}{\tp[S]{X}{N}}{\Hom[S]{X^*}{N}},
  \end{equation*}
  given by evaluation $\theta(x\ot n)(\z) = \z(x)n$. Auslander
  computed the kernel and cokernel of this map in
  \prpcite[6.3]{MAs66}. Because the map is pivotal for our proof of
  the next lemma, we include a computation for the case where $X$ is
  totally reflexive.
  
  Consider a short exact sequence $0 \to N_1 \to F \to N \to 0$, where
  $F$ is a free $S$-module. For any totally reflexive $S$-module $X$,
  the evaluation homomorphism $\theta_{XF}$ is an isomorphism, and the
  commutative diagram
  \begin{equation*}
    \xymatrix@C=1em{
      & \tp[S]{X}{\syz{1}{N}}
      \ar[r]
      \ar[d]^-{\theta_{X\syz{1}{N}}}
      & \tp[S]{X}{F}
      \ar[r]
      \ar[d]_-{\is}^-{\theta_{XF}}
      & \tp[S]{X}{N}
      \ar[r]
      \ar[d]^-{\theta_{XN}}
      & 0
      \\
      0
      \ar[r]
      & \Hom[S]{X^*}{\syz{1}{N}}
      \ar[r]
      & \Hom[S]{X^*}{F}
      \ar[r]
      & \Hom[S]{X^*}{N}
      \ar[r]
      & \Ext[S]{1}{X^*}{\syz{1}{N}} \ar[r] & 0
    }
  \end{equation*}
  shows that there is an isomorphism $\Coker{\theta_{XN}} \is
  \Ext[S]{1}{X^*}{\syz{1}{N}}$.  The snake lemma applies to yield
  $\Ker{\theta_{XN}} \is \Coker{\theta_{X\syz{1}{N}}} \is
  \Ext[S]{1}{X^*}{\syz{2}{N}}$, and then \eqref{dim shift} gives
  \begin{equation}
    \label{eq:Ker-Coker}
    \Ker{\theta_{XN}} \is \Ext[S]{1}{\syzp{2}{X}^*}{N} \qand 
    \Coker{\theta_{XN}} \is \Ext[S]{1}{\syzp{1}{X}^*}{N}.
  \end{equation}
\end{ipg}

\begin{lem}
  \label{lem:shopping list}
  Let $\Snl$ be a complete local ring of depth $0$.  Let $\cat{C}$ be
  a reflexive subcategory of $\,\mod{S}$.  If $\l$ has a
  $\cat{C}$-approximation and $\l$ is not in $\cat{C}$, then $\cat{C}
  = \F[S]$.
\end{lem}

\begin{prf*}
  Consider a $\cat{C}$-approximation $0 \lra L \xra{\a} C \lra \l \lra
  0$, and dualize to get $0 \lra \l^* \lra C^* \xra{\a^*} L^*$.  Let
  $I$ be the image of $\a^*$, and let $\f$ be the factorization of
  $\a^*$ through the inclusion $I \into L^*$.
  
  First we prove that the surjection $\f$ is a $\cl{\cat{C}}$-precover
  of $I$. Let $X$ be a module in $\cl{\cat{C}}$. If $X$ is a free
  $S$-module, then any homomorphism $X \to I$ lifts through $\f$. We
  may now assume that $X$ is indecomposable and not free. Because
  $\Hom[S]{X}{I}$ is a submodule of $\Hom[S]{X}{L^*}$, it suffices to
  prove surjectivity of
  \begin{equation*}
    \dmapdef{\Hom[S]{X}{\a^*}}{\Hom[S]{X}{C^*}}{\Hom[S]{X}{L^*}},
  \end{equation*}
  which we do next.

  The vertical maps in the commutative diagram below are evaluation
  homomorphisms, see~\pgref{nat}.
  \begin{equation*}
    \xymatrix@C=1.5em{
      & \tp{X}{L} 
      \ar[r]^-\iota
      \ar[d]_-{\theta_{XL}}
      & \tp{X}{C} 
      \ar[r] 
      \ar[d]_-{\theta_{XC}}
      & \tp{X}{\l} 
      \ar[r] 
      \ar[d]_-{\theta_{X\l}}
      & 0
      \\
      0 
      \ar[r]       
      & \Hom[S]{X^*}{L} 
      \ar[r] 
      & \Hom[S]{X^*}{C} 
      \ar[r]
      & \Hom[S]{X^*}{\l}
      \ar[r]
      & \Ext[S]{1}{X^*}{L}
    }
  \end{equation*}
  First we argue that the rows of this diagram are short exact
  sequences.  The module $X$ is in $\cl{\cat{C}}$ and hence in
  $\G[S]$, see~\pgref{FBG}, so $\Ext[S]{1}{X^*}{L}=0$. Moreover,
  $\theta_{XL}$ is an isomorphism by \eqref{Ker-Coker}, hence $\iota$
  is injective. Next note that for every $\z \in X^*$ the image of
  $\mapdef{\z}{X}{S}$ is in $\n$ as $X$ is indecomposable and not
  free. Thus, for all $x\in X$ and $u\in \l$, we have
  $\theta_{X\l}(x\ot u)(\z) = \z(x)u=0$. Finally, apply
  $\Hom[S]{-}{S}$ to the diagram above and use Hom-tensor adjointness
  to get
  \begin{equation*}
    \xymatrix@C=1.5em{
      \Hom[S]{X^*}{\l}^* 
      \ar[r] 
      \ar[d]_-{0}
      & \Hom[S]{X^*}{C}^* 
      \ar[r]
      \ar[d]_-{\theta_{XC}^*}
      & \Hom[S]{X^*}{L}^*
      \ar[r]
      \ar[d]^-{\is}_-{\theta_{XL}^*}
      & \Ext[S]{1}{\Hom[S]{X^*}{\l}}{S}
      \ar[d]_-{0}
      \\
      \Hom[S]{X}{\l^*}
      \ar[r]
      & \Hom[S]{X}{C^*} 
      \ar[r]^-{\stackrel{\stackrel{\Hom[S]{X}{\a^*}}{}}{}} 
      & \Hom[S]{X}{L^*} 
      \ar[r] 
      & \Ext[S]{1}{\tp{X}{\l}}{S}.
    }
  \end{equation*}
  The diagram shows that $\Hom[S]{X}{\a^*}$ is surjective, as desired.
  
  Now $\mapdef{\f}{C^*}{I}$ is a $\cl{\cat{C}}$-precover, so by
  completeness of $S$, the module $I$ has a $\cl{\cat{C}}$-cover; see
  \pgref{cov-approx}(c).  The ring has depth $0$, so $\l^*$ is a
  non-zero $\l$-vector space. By the assumptions on $\cat{C}$, the
  residue field $\l$ cannot be a direct summand of $C^*$. As $\Ker{\f}
  = \l^*$, it follows from \lemref{xu} that $\f$ is a
  $\cl{\cat{C}}$-cover.  For every $X\in \cl{C}$ Wakamatsu's lemma
  gives $\Ext[S]{1}{X}{\l^*}=0$. Consequently, every module in
  $\cat{C}$ is projective and hence free, since $S$ is local.
\end{prf*}

\begin{ipg}
  \label{GMC}
  Let $\Rmk$ be a local ring and denote by $\cat{M}(R)$ the category
  of maximal Cohen--Macaulay $R$-modules.
  
  (a) If $R$ is Cohen--Macaulay, then $\G \subseteq \cat{M}(R)$ by the
  Auslander--Bridger formula \cite[\S 3.2 Prop.~3]{MAs67}.
  Conversely, if $\G \subseteq \cat{M}(R)$, then $R$ is
  Cohen--Macaulay.
  
  (b) If $R$ is Gorenstein, then the categories $\G$ and $\cat{M}(R)$
  coincide by \cite[\S 3.2 Thm.~3]{MAs67} and the Auslander--Bridger
  formula.  Conversely, if $\G = \cat{M}(R)$, then $R$ is Gorenstein.
  Indeed, $R$ is Cohen--Macaulay by (a), so $\Syz{\dimR}{\k}$ is in
  $\cat{M}(R)$, hence in $\G$, and therefore $R$ is Gorenstein by
  \cite[\S 3.2, Rmk.\ after Thm.\ 3]{MAs67}.
  
  (c) If $R$ is Gorenstein, then a short exact sequence $0 \to L \to G
  \to M \to 0$ is a CM-approximation if and only if it is a
  $\G$-approximation. This follows from (b) and the fact that $L$ is
  in $\cat{M}(R)^\perp$ if and only if $L$ has finite injective
  dimension.
\end{ipg}

If $R$ is Gorenstein, then every $R$-module has a CM-approximation by
\thmcite[A]{MAsROB89}. In view of \pgref{GMC}(c) the next result
contains a converse, cf.~Theorem C.

\begin{thm}
  \label{thm:B-approximation}
  Let $\Rmk$ be a local ring and $\cat{B}$ be a reflexive subcategory
  of $\,\mod{R}$. If $\,\k$ has a $\cat{B}$-approximation, then $R$ is
  Gorenstein or $\cat{B} = \F$.
\end{thm}

In our proof of this theorem we use the next lemma. We do not know a
reference giving a direct argument, so one is supplied here.

\begin{lem}
  \label{lem:k}
  Let $\Rmk$ be a local ring, and let $\,\x = x_1,\dots,x_n$ be a
  sequence in $\m \setminus \m^2$. If $\x$ is linearly independent
  modulo $\m^2$, then $\k$ is a direct summand of the module
  $\Syz{n}{\k}/\x\Syz{n}{\k}$.
\end{lem}

\begin{prf*}
  Let $(\Kc(\x),\ddK)$ be the Koszul complex on $\x$. If necessary,
  supplement $\x$ to a minimal generating sequence $\x,\y$ for $\m$.
  Let $(\pmb{F},\dd)$ be a minimal free resolution of~$\k$.  The
  identification $R/(\x,\y) = \k$ lifts to a morphism of complexes
  $\mapdef{\sigma}{\Kc(\x,\y)}{\pmb{F}}$. Serre proves in
  \cite[Appendix~I.2]{localg} that $\sigma$ is injective and
  degreewise split. The natural inclusion
  $\mapdef[\into]{\iota}{\Kc(\x)}{\Kc(\x,\y)}$ is also degreewise
  split, so the composite $\rho = \sigma\iota$ is an injective
  morphism of complexes and degreewise split.
  
  From the short exact sequence $0 \lra \Syz{n}{\k} \xra{\iota}
  F_{n-1} \lra \Syz{n-1}{\k} \lra 0$, we get an exact sequence in
  homology that reads in part
  \begin{equation*}
    \tag{\one}
    \Tor{1}{R/(\x)}{\Syz{n-1}{\k}} \lra \tp{R/(\x)}{\Syz{n}{\k}}
    \xra{\smash{\tp{R/(\x)}{\iota}}} \tp{R/(\x)}{F_{n-1}}.
  \end{equation*}
  The module $\Tor{1}{R/(\x)}{\Syz{n-1}{\k}} \is \Tor{n}{R/(\x)}{\k}$
  is annihilated by $\m$.
  
  Let $e$ be a generator of $\Kc(\x)_n$. The image $\rho_n(e)$ in
  $F_n$ is a minimal generator as $\rho_n$ is split.  Set $\be =
  \dd_n\rho_n(e) \in \Syz{n}{\k}$; since $\pmb{F}$ is minimal, $\be$
  is a minimal generator of the syzygy $\Syz{n}{\k}$.  The minimal
  generator $\tp[]{1}{\be}$ of $\tp{R/(\x)}{\Syz{n}{\k}}$ is in the
  kernel of $\tpp{R/(\x)}{\iota}$, as the element $\be =
  \dd_n\rho_n(e) = \rho_{n-1}\ddK_n(e)$ is in $\x F_{n-1}$.  By
  exactness of (\one) the element $\tp[]{1}{\be}$ is annihilated by
  $\m$, hence it generates a \mbox{$1$-dimensional} $\k$-vector space
  that is a direct summand of $\Syz{n}{\k}/\x\Syz{n}{\k}$.
\end{prf*}

\begin{prf*}[Proof of \pgref{thm:B-approximation}]
  We aim to apply \lemref{shopping list}. By
  \prpref[Propositions~]{C-cover} and \prpref[]{laundry list}, and by
  faithful flatness of $\Rhat$, we may assume $R$ is complete. Set
  $d=\dptR$; by \lemref{syz-approximation} the $d$th syzygy
  $\Syz{d}{\k}$ has a $\cat{B}$-approximation:
  \begin{equation*}
    0 \to L \to B \to \Syz{d}{\k} \to 0.
  \end{equation*}
  Let $\x = x_1,\dots,x_d$ be an $R$-regular sequence in $\m \setminus
  \m^2$ linearly independent modulo $\m^2$. The Koszul homology
  modules
  $$\operatorname{H}_i(\tp{\Kc(\x)}{\Syz{d}{\k}}) \is
  \Tor{i}{R/(\x)}{\Syz{d}{\k}} \is \Tor{i+d}{R/(\x)}{\k}$$ vanish for
  $i>0$, so $\x$ is also $\Syz{d}{\k}$-regular.
  
  Set $S = R/(\x)$; by \pgref{FBG} and \prpref{C-cover} the sequence
  \begin{equation*}
    0 \lra \tp{S}{L} \lra \tp{S}{B} \xra{\varphi}
    \tp{S}{\Syz{d}{\k}} \lra 0
  \end{equation*}
  is a $\cl{\bar{B}}$-approximation. Moreover, the category
  $\cl{\bar{B}}$ is reflexive by \prpref{laundry list}. By \lemref{k}
  the residue field $\k$ is a direct summand of $\tp{S}{\Syz
    {d}{\k}}$, so by \pgref{summand} there is an
  $\cl{\bar{B}}$-precover of $\k$. Since $S$ is complete, it follows
  from \pgref{precov-approx} that $\k$ has a
  $\cl{\bar{B}}$-approximation.
  
  Assume $R$ is not Gorenstein. Then $S$ is not Gorenstein, so the
  residue field $\k$ is not in $\G[S]$ and hence not in
  $\cl{\bar{B}}$; see \cite[\S 3.2, Rmk.\ after Thm.\ 3]{MAs67} or
  \thmcite[(1.4.9)]{lnm}. By \lemref{shopping list} every module in
  $\cl{\bar{B}}$ is now free, so for every $B\in\cat{B}$ the module
  $\tp{S}{B}$ is free over $S$.  By \pgref{FBG} the sequence $\x$ is
  $B$-regular; therefore, $B$ is a free $R$-module by Nakayama's
  lemma.
\end{prf*}

An approximation of a module $M$ is \emph{minimal} if the map onto $M$
is a cover. When $R$ is Gorenstein, every $R$-module has a minimal
CM-approximation by unpublished work of Auslander; see
\cite[Sec.~4]{ADS-93} and \thmcite[5.5]{EJX-99}. Hence we have

\begin{cor}
  \label{cor:G-cover}
  Let $\Rmk$ be a local ring and assume there is a non-free module in
  $\G$. The following are then equivalent:
  \begin{eqc}
  \item $R$ is Gorenstein.
  \item $\k$ has a $\G$-approximation.
  \item Every finitely generated $R$-module has a minimal
    $\G$-approximation. \qed
  \end{eqc}
\end{cor}

\begin{ipg}
  If $R$ has a dualizing complex, cf.~\cite[V.\S 2]{rad}, then $\k$
  has a Gorenstein projective precover $X \to \k$ by
  \thmcite[2.11]{PJr07}. Assume $X$ is finitely generated, i.e., $X$
  is in $\G$ and, further, that $R$ is henselian. If $X$ is free, then
  it follows from \pgref{precov-approx} that $\k$ has a
  $\G$-approximation $0 \to L \to X' \to \k \to 0$, where $X'$ is
  free. Hence, $\k$ is in $\G^\perp$ and then $\G = \F$.  If $X$ is
  not free, then $R$ is Gorenstein by \corref[]{G-cover}.
\end{ipg}

\begin{bfhpg}[Questions]
  \label{q1}
  Let $\Rmk$ be a local ring. If $\k$ has a $\G$-precover, is then
  $\G$ precovering?  If $\G$ is precovering and contains a non-free
  module, is then $R$ Gorenstein?
\end{bfhpg}

\section{On the number of totally reflexive modules}
\label{sec:number}

In this section we prove Theorems~A and B. Note that by \pgref{pc} the
latter would follow immediately from a positive answer to the second
question in \pgref{q1}.

\begin{lem}
  \label{lem:finite length}
  Let $R$ be a local ring and $M$ and $N$ be finitely generated
  $R$-modules. If only finitely many isomorphism classes of
  $R$-modules $X$ can fit in a short exact sequence $0\to N \to X\to M
  \to 0$, then the $R$-module $\Ext{1}{M}{N}$ has finite length.
\end{lem}

\begin{prf*}
  Given an $R$-module $X$, we denote by $[X]$ the subset of
  $\Ext{1}{M}{N}$ whose elements have representatives of the form $0
  \to N \to Y \to M \to 0$, where \mbox{$Y \cong X$}.  By assumption,
  there exist non-isomorphic $R$-modules $X_0, \dots, X_n$ such
  that\linebreak[4] $\Ext{1}{M}{N}$ is the disjoint union of the sets
  $[X_i]$.  We may take $X_0=M\oplus N$, so $[X_0]$ is the zero
  submodule of $\Ext{1}{M}{N}$.  We must prove that there is an
  integer $q>0$ such that $\m^q\Ext{1}{M}{N}$ is contained in $[X_0]$.
  
  By \corcite[1]{RMG85} there are integers $p_i$ such that if
  $M/\m^pM\oplus N/\m^pN \is X_i/\m^pX_i$ for some $p \ge p_i$, then
  $X_i\is M\oplus N$. Set $q = \max\set{p_1,\dots,p_n}$. Take a short
  exact sequence $\xi$ in $\m^q\Ext{1}{M}{N}$; it belongs to some set
  $[X_i]$. By \thmcite[1.1]{JSt05} the sequence $\tp{\xi}{R/\m^q}$
  splits, so $M/\m^qM\oplus N/\m^qN \is X_i/\m^qX_i$. By the choice of
  $q$ this implies $X_i\is M\oplus N$, so $i=0$, i.e.\ $\xi$ is in the
  zero submodule $[X_0]$.
\end{prf*}

Let $R\to S$ be a flat ring homomorphism. It does not follow from the
natural isomorphism $\tp{S}{\Ext{1}{M}{N}} \is
\Ext[S]{1}{\tp{S}{M}}{\tp{S}{N}}$ that every extension of the
$S$-modules $\tp{S}{N}$ and $\tp{S}{M}$ has the form $\tp{S}{X}$ for
some $R$-module $X$. In a seminar, Roger Wiegand alerted us to the
next result.

\begin{lem}
  \label{lem:from roger}
  Let $(R,\m) \to (S,\n)$ be a flat ring homomorphism with $\m S=\n$
  and $R/\m \is S/\n$.  Let $M$ and $N$ be finitely generated
  $R$-modules and $\xi$ be an element of the $S$-module
  $\Ext[S]{1}{\tp{S}{M}}{\tp{S}{N}}$. If the $R$-module
  $\Ext{1}{M}{N}$ has finite length, then there is an element $\chi$
  in $\Ext{1}{M}{N}$ such that $\xi = \tp{S}{\chi}$.
\end{lem}

\begin{prf*}
  The functor $\tp{S}{-}$ from the category $\mod{R}$ to itself
  induces a natural isomorphism $K \to \tp{S}{K}$ on $R$-modules of
  finite length. Applied to $\Ext{1}{M}{N}$ this yields the first
  isomorphism below
  \begin{equation*}
    \Ext{1}{M}{N} \ira \tp{S}{\Ext{1}{M}{N}} \ira
    \Ext[S]{1}{\tp{S}{M}}{\tp{S}{N}}.
  \end{equation*}
  The composite sends an exact sequence $\chi$ to $\tp{S}{\chi}$.
\end{prf*}

The next result is Theorem~B from the introduction.

\begin{thm}
  \label{thm:B}
  Let $R$ be a local ring. If the set of isomorphism classes of
  indecomposable modules in $\G$ is finite, then $R$ is Gorenstein or
  $\G = \F$.
\end{thm}

\begin{prf*}
  Assume there are only finitely many isomorphism classes of
  indecomposable modules in $\G$.  By \pgref{pc} the residue field
  $\k$ then has a $\G$-precover $\mapdef{\f}{B}{\k}$. We claim that
  $\tp{\Rhat}{\f}$ is a $\cl{\tp[]{\Rhat}{\G}}$-precover of $\k$.
  Since $\Rhat$ is complete, this implies the existence of a
  $\cl{\tp[]{\Rhat}{\G}}$-approximation of $\k$, see
  \pgref{precov-approx}, and the desired conclusion follows from
  \thmref{B-approximation} and faithful flatness of $\Rhat$.
  
  To prove the claim, we must show that
  \begin{equation*}
    \dmapdef{\Hom[\Rhat]{H'}{\tp{\Rhat}{\f}}}{\Hom[\Rhat]{H'}{\tp{\Rhat}{B}}}
    {\Hom[\Rhat]{H'}{\k}}
  \end{equation*}
  is surjective for every module $H'\in \cl{\tp[]{\Rhat}{\G}}$. By
  flatness of $\Rhat$, surjectivity holds for modules in
  $\tp[]{\Rhat}{\G}$ and hence for every module in
  $\add{\tp[]{\Rhat}{\G}}$. It is now sufficient to prove that the
  category $\add{\tp[]{\Rhat}{\G}}$ is closed under extensions,
  because then $\cl{\tp[]{\Rhat}{\G}}$ is $\add{\tp[]{\Rhat}{\G}}$.
  
  First we show that $\tp[]{\Rhat}{\G}$ is closed under extensions.
  Fix modules $G$ and $K$ in $\G$, and consider short exact sequences
  $0 \to G \to H \to K \to 0$. Each $H$ is in $\G$, and the minimal
  number of generators of each $H$ is bounded by the sum of the
  numbers of minimal generators for $G$ and $K$. Since the number of
  indecomposable modules in $\G$ is finite, there are, up to
  isomorphism, only finitely many such modules $H$. By \lemref{finite
    length} the module $\Ext{1}{K}{G}$ has finite length, and by
  \lemref[]{from roger} every element of
  $\Ext[\Rhat]{1}{\tp{\Rhat}{K}}{\tp{\Rhat}{G}}$ is extended from
  $\Ext{1}{K}{G}$.

  To prove that $\add{\tp[]{\Rhat}{\G}}$ is closed under extensions,
  let $G'$ and $K'$ be summands of extended modules, i.e.,\ $G'\oplus
  G'' \is \tp{\Rhat}{G}$ and $K'\oplus K''\is \tp{\Rhat}{K}$ for
  modules $G,K \in\G$. Consider a short exact sequence $0 \to G' \to
  H' \to K' \to 0$.  Then a sequence
  \begin{equation*}
    0 \lra G'\oplus G'' \lra H'\oplus G''\oplus K'' \lra K'\oplus K''
    \lra 0,
  \end{equation*}
  is exact, so by what has already been proved, the middle term
  $H'\oplus G''\oplus K''$ is in $\tp[]{\Rhat}{\G}$; whence $H'$ is in
  $\add{\tp[]{\Rhat}{\G}}$.
\end{prf*}

In view of \pgref{GMC}(a) we have

\begin{cor}
  \label{cor:finiteCMtype}
  Let $R$ be a Cohen--Macaulay local ring. If $R$ is of finite CM
  representation type, then $R$ is Gorenstein or $\G = \F$. \qed
\end{cor}

The next result contains Theorem~A from the introduction.

\begin{thm}
  \label{thm:A}
  Let $R$ be a local ring and assume the set of isomorphism classes of
  indecomposable modules in $\G \setminus \F$ is finite and not empty.
  Then $R$ is Gorenstein and an isolated singularity.  Further,
  $\Rhat$ is a hypersurface singularity; if finite CM representation
  type ascends from $R$ to $\Rhat$, then $\Rhat$ is even a
  simple~singularity.
\end{thm}

\begin{prf*}
  By \thmref{B} the ring $R$ is Gorenstein. From \pgref{GMC}(b) it
  follows that $R$ is of finite CM representation type and hence an
  isolated singularity by \corcite[2]{CHnGJL02}.  By \cite[Satz
  1.2]{JHr78} the completion $\Rhat$ is a hypersurface singularity
  and, assuming that also $\Rhat$ is of finite CM representation type,
  it follows from \corcite[(8.16)]{yos} that $\Rhat$ is a \mbox{simple
    singularity.}
\end{prf*}

\begin{rmk}
  In \cite{FOS87} Schreyer conjectured that a Cohen--Macaulay local
  $\k$-algebra $R$ is of finite CM representation type if and only if
  $\Rhat$ is of finite CM representation type. In \cite{RWg98}
  R.~Wiegand proved descent of finite CM representation type from
  $\Rhat$ to $R$ for any local ring $R$. Ascent is verified in
  \cite{RWg98} when $R$ is Cohen--Macaulay and either $\Rhat$ is an
  isolated singularity or $\dimR \le 1$. Ascent also holds for
  excellent Cohen--Macaulay local rings by work of Leuschke and
  R.~Wiegand \cite{GLsRWg00}.
\end{rmk}

\begin{bfhpg}[Remarks]
  Constructing rings with infinitely many totally reflexive modules is
  easy using \thmref{B}. Indeed, let $Q$ be a local ring of positive
  dimension and set $R=\pows[Q]{X}/(X^2)$. As $R$ is not reduced, it
  is not an isolated singularity.  The $R$-module $R/(X)$ is in $\G$
  and is not free, cf.~\cite[exa.~(4.1.5)]{lnm}, so by \thmref[]{B}
  there are infinitely many non-isomorphic indecomposable modules in
  $\G$.
  
  More generally, Avramov, Gasharov, and Peeva \cite{AGP-97} construct
  a non-free totally reflexive module\footnote{ Actually, even a
    module of CI-dimension $0$ as defined in \cite[(1.2)]{AGP-97}.}
  $G$ over any ring of the form $R \is Q/(\x)$, where $(Q,\q)$ is
  local and $\x \in \q^2$ is a $Q$-regular sequence. Such a ring $R$
  is said to have an embedded deformation of codimension $c$, where
  $c$ is the length of $\x$. Again \thmref[]{B} implies the existence
  of infinitely many non-isomorphic indecomposable modules in $\G$. If
  $\Rhat$ has an embedded deformation of codimension $c \ge 2$, a
  recent argument of Avramov and Iyengar builds from $G$ an infinite
  family of non-isomorphic indecomposable modules in $\G$; see
  \thmcite[6.8 and proof of 6.4.(1)]{LLASIn07}. For such $R$, this gives
  a constructive proof of the abundance of modules in $\G$.
\end{bfhpg}

\begin{bfhpg}[Question]
  Let $R$ be a local ring that is not Gorenstein. Given an
  indecomposable totally reflexive $R$-module $G\not\is R$, are there
  constructions that produce infinite families of non-isomorphic
  indecomposable modules in $\G$?
\end{bfhpg}


\section*{Acknowledgments}

We thank Luchezar Avramov for supporting R.T.'s visit to the
University of Nebraska and for advice in the process of writing up
this work.



\def\cprime{$'$}
  \newcommand{\arxiv}[2][AC]{\mbox{\href{http://arxiv.org/abs/#2}{\sf arXiv:#2
  [math.#1]}}}
  \newcommand{\oldarxiv}[2][AC]{\mbox{\href{http://arxiv.org/abs/math/#2}{\sf
  arXiv:math/#2
  [math.#1]}}}\providecommand{\MR}[1]{\mbox{\href{http://www.ams.org/mathscine%
t-getitem?mr=#1}{MR#1}}}
  \renewcommand{\MR}[1]{\mbox{\href{http://www.ams.org/mathscinet-getitem?mr=#%
1}{MR#1}}}
\providecommand{\bysame}{\leavevmode\hbox to3em{\hrulefill}\thinspace}
\providecommand{\MR}{\relax\ifhmode\unskip\space\fi MR }
\providecommand{\MRhref}[2]{%
  \href{http://www.ams.org/mathscinet-getitem?mr=#1}{#2}
}
\providecommand{\href}[2]{#2}

\end{document}